\documentclass[10pt]{article}
\usepackage{amsmath}
\usepackage{latexsym}
\usepackage{amssymb}
\usepackage{amsfonts}
\usepackage{amsthm}
\title{\bf{THE MATHEMATICAL FOUNDATIONS OF \\  GAUGE THEORY REVISITED}}
\author{\bf{Jean-Francois Pommaret }\\ CERMICS, Ecole des Ponts ParisTech, 77455 Marne-la-Vallee, France \\
Email: jean-francois.pommaret@wanadoo.fr, pommaret@cermics.enpc.fr   }
\date{  }
\begin{document}
\maketitle

\noindent
{\bf ABSTRACT}  \\

We start recalling with critical eyes the mathematical methods used in gauge theory and prove that they are not coherent with continuum mechanics, in particular the analytical mechanics of rigid bodies or hydrodynamics, though using the same group theoretical methods and despite the well known couplings existing between elasticity and electromagnetism (piezzoelectricity, photoelasticity, streaming birefringence). The purpose of this paper is to avoid such contradictions by using new mathematical methods coming from the formal theory of systems of partial differential equations and Lie pseudogroups. These results finally allow to unify the previous independent tentatives done by the brothers E. and F. Cosserat in 1909 for elasticity or H. Weyl in 1918 for electromagnetism by using respectively the group of rigid motions of space or the conformal group of space-time. Meanwhile we explain why the Poincar\'{e} {\it duality scheme} existing between {\it geometry} and {\it physics} has to do with homological algebra and algebraic analysis. We insist on the fact that these results could not have been obtained before 1975 as the corresponding tools were not known before.\\

\noindent
{\bf Keywords}: Gauge Theory, Curvature, Torsion, Maurer-Cartan Forms, Maurer-Cartan Equations, Lie Groups, Lie Pseudogroups, Differential Sequence, Poincar\'{e} sequence, Janet Sequence, Spencer Sequence, Differential Modules, Homological Algebra, Extension Modules. \\

\noindent
{\bf 1. Introduction} \\

It is usually accepted today in the literature that the {\it physical foundations} of what we shall simply call (classical) "gauge theory" (GT) can be found in the paper [31] published by C.N. Yang and R.L. Mills in 1954. Accordingly, the {\it mathematical foundations} of GT can be found in the references existing at this time on differential geometry and group theory, the best and most quoted one being the survey book [9] published by S. Kobayashi and K. Nomizu in 1963 (See also [4,6,7,30]). The pupose of this introduction is to revisit these foundations with critical eyes, recalling them in a quite specific and self-contained way for later purposes. \\

The word "{\it group}" has been introduced for the first time in 1830 by Evariste Galois (1811-1832). Then this concept slowly passed from algebra (groups of permutations) to geometry (groups of transformations). It is only in 1880 that Sophus Lie (1842-1899) studied the groups of transformations depending on a finite number of parameters and now called {\it Lie groups of transformations}. \\

 Let $X$ be a manifold with local coordinates $x=(x^1, ... , x^n)$ and $G$ be a {\it Lie group}, that is another manifold with local coordinates $a=(a^1, ... , a^p)$ called {\it parameters} with a {\it composition} $G\times G \rightarrow G: (a,b)\rightarrow ab$, an {\it inverse} $G \rightarrow G: a \rightarrow a^{-1}$ and an {\it identity} $e\in G$ satisfying:\\
\[(ab)c=a(bc)=abc,\hspace{1cm} aa^{-1}=a^{-1}a=e,\hspace{1cm} ae=ea=a,\hspace{1cm} \forall a,b,c \in G \]
Then $G$ is said to {\it act} on $X$ if there is a map $X\times G \rightarrow X: (x,a) \rightarrow y=ax=f(x,a)$ such that $(ab)x=a(bx)=abx, \forall a,b\in G, \forall x\in X$ and, for simplifying the notations, we shall use global notations even if only local actions are existing. The action is said to be {\it effective} if $ax=x, \forall x\in X\Rightarrow a=e$. A subset $S\subset X$ is said to be {\it invariant} under the action of $G$ if $aS\subset S,\forall a\in G$ and the {\it orbit} of $x\in X$ is the invariant subset $Gx=\{ax\mid a\in G\}\subset X$. If $G$ acts on two manifolds $X$ and $Y$, a map $f:X\rightarrow Y$ is said to be {\it equivariant} if $f(ax)=af(x), \forall x\in X, \forall a\in G$. For reasons that will become clear later on, it is often convenient to introduce the {\it graph} $X\times G\rightarrow X\times X: (x,a)\rightarrow (x,y=ax)$ of the action. In the product $X\times X$, the first factor is called the {\it source} while the second factor is called the {\it target}. \\

We denote as usual by $T=T(X)$ the {\it tangent bundle} of $X$, by $T^*=T^*(X)$ the {\it cotangent bundle}, by ${\wedge}^rT^*$ the {\it bundle of r-forms} and by $S_qT^*$ the {\it bundle of q-symmetric tensors}. Moreover, if  $\xi,\eta\in T$ are two vector fields on $X$, we may define their {\it bracket} $[\xi,\eta]\in T$ by the local formula $([\xi,\eta])^i(x)={\xi}^r(x){\partial}_r{\eta}^i(x)-{\eta}^s(x){\partial}_s{\xi}^i(x)$ leading to the {\it Jacobi identity} $[\xi,[\eta,\zeta]]+[\eta,[\zeta,\xi]]+[\zeta,[\xi,\eta]]=0, \forall \xi,\eta,\zeta \in T$ allowing to define a {\it Lie algebra}. We have also the useful formula $[T(f)(\xi),T(f)(\eta)]=T(f)([\xi,\eta])$ where $T(f):T(X)\rightarrow T(Y)$ is the tangent mapping of a map $f:X\rightarrow Y$. Finally, when $I=\{ i_1< ... < i_r\}$ is a multi-index, we may set $dx^I=dx^{i_1}\wedge ... \wedge dx^{i_r}$ and introduce the {\it exterior derivative} $d:{\wedge}^rT^*\rightarrow {\wedge}^{r+1}T^*:\omega={\omega}_Idx^I \rightarrow d\omega={\partial}_i{\omega}_Idx^i\wedge dx^I$ with $d^2=d\circ d\equiv 0$ in the {\it Poincar\'{e} sequence}:\\
\[  {\wedge}^0T^* \stackrel{d}{\longrightarrow} {\wedge}^1T^* \stackrel{d}{\longrightarrow} {\wedge}^2T^* \stackrel{d}{\longrightarrow} ... \stackrel{d}{\longrightarrow} {\wedge}^nT^* \longrightarrow 0  \]

In order to fix the notations, we quote without any proof the "{\it three fundamental theorems of Lie}" that will be of constant use in the sequel (See [17] for more details):\\

\noindent
{\bf FIRST FUNDAMENTAL THEOREM  1.1}: The orbits $x=f(x_0,a)$ satisfy the system of PD equations $\partial x^i/\partial a^{\sigma}= {\theta}^i_{\rho}(x){\omega}^{\rho}_{\sigma}(a)$  with $det(\omega)\neq 0$. The vector fields ${\theta}_{\rho}={\theta}^i_{\rho}(x){\partial}_i$ are called {\it infinitesimal generators} of the action and are linearly independent over the constants when the action is effective.\\

In a rough way, we have $x=ax_0\Rightarrow dx=dax_0=daa^{-1}x$ and $daa^{-1}=\omega=({\omega}^{\tau}={\omega}^{\tau}_{\sigma}(a)da^{\sigma})$ is thus a family of right invariant 1-forms on $G$ called {\it Maurer-Cartan forms} or simply MC {\it forms}.  \\
 
\noindent
{\bf SECOND FUNDAMENTAL THEOREM  1.2}: If ${\theta}_1,...,{\theta}_p$ are the infinitesimal generators of the effective action of a lie group $G$ on $X$, then $[{\theta}_{\rho},{\theta}_{\sigma}]=c^{\tau}_{\rho\sigma}{\theta}_{\tau}$ where the $c=(c^{\tau}_{\rho\sigma}= - c^{\tau}_{\sigma\rho})$ are the {\it structure constants} of a Lie algebra of vector fields which can be identified with ${\cal{G}}=T_e(G)$ the tangent space to $\cal{G}$ at the identity $e\in G$ by using the action as we already did. Equivalently, introducing the non-degenerate inverse matrix $\alpha={\omega}^{-1}$ of right invariant vector fields on $G$, we obtain from crossed-derivatives the {\it compatibility conditions} (CC) for the previous system of partial differential (PD) equations called {\it Maurer-Cartan equations} or simply MC {\it equations}, namely: \\
\[  \partial {\omega}^{\tau}_s/\partial a^r - \partial {\omega}^{\tau}_r/ \partial a^s + c^{\tau}_{\rho\sigma} {\omega}^{\rho}_r {\omega}^{\sigma}_s  = 0\]
({\it care to the sign used}) or equivalently $[{\alpha}_{\rho},{\alpha}_{\sigma}]=c^{\tau}_{\rho\sigma} {\alpha}_{\tau} $.  \\
 
 Using again crossed-derivatives, we obtain the corresponding {\it integrability conditions} (IC) on the structure constants and the Cauchy-Kowaleski 
 theorem finally provides:\\
 
\noindent 
{\bf THIRD FUNDAMENTAL THEOREM  1.3}: For any Lie algebra $\cal{G}$ defined by structure constants $c=(c^{\tau}_{\rho\sigma})$ satisfying :\\
\[  c^{\tau}_{\rho\sigma}+c^{\tau}_{\sigma \rho}=0, \hspace{1cm} c^{\lambda}_{\mu\rho}c^{\mu}_{\sigma\tau}+c^{\lambda}_{\mu\sigma}c^{\mu}_{\tau\rho}+c^{\lambda}_{\mu\tau}c^{\mu}_{\rho\sigma}=0  \]
one can construct an analytic group $G$ such that ${\cal{G}}=T_e(G)$ by recovering the MC forms from the MC equations.\\

\noindent
{\bf EXAMPLE  1.4}: Considering the affine group of transformations of the real line $y=a^1x+a^2$, the orbits are defined by $x=a^1x_0+a^2$, a definition leading to $dx=((1/a^1)da^1)x+(da^2-(a^2/a^1)da^1)$. We obtain therefore ${\theta}_1=x{\partial}_x, {\theta}_2={\partial}_x \Rightarrow [{\theta}_1,{\theta}_2]=-{\theta}_2$ and ${\omega}^1=(1/{a^1})da^1, {\omega}^2=da^2-(a^2/{a^1})da^1\Rightarrow d{\omega}^1=0, d{\omega}^2-{\omega}^1\wedge{\omega}^2=0 \Leftrightarrow [{\alpha}_1,{\alpha}_2]=-{\alpha}_2$ with ${\alpha}_1=a^1{\partial}_1+a^2{\partial}_2, {\alpha}_2={\partial}_2$.\\
 
\noindent
{\bf GAUGING PROCEDURE  1.5}: If $x=a(t)x_0+b(t)$ with $a(t)$ a time depending orthogonal matrix ({\it rotation}) and $b(t)$ a time depending vector ({\it translation}) describes the movement of a rigid body in ${\mathbb{R}}^3$, then the projection of the {\it absolute speed} $v=\dot{a}(t)x_0+\dot{b}(t)$ in an orthogonal frame fixed in the body is the so-called {\it relative speed} $a^{-1}v=a^{-1}\dot{a}x_0+a^{-1}\dot{b}$ and the kinetic energy/Lagrangian is a quadratic function of the $1$-forms $A=(a^{-1}\dot{a}$, $a^{-1}\dot{b})$. Meanwhile, taking into account the preceding example, the {\it Eulerian speed}  $v=v(x,t)=\do{a}a^{-1}x+\dot{b}-\dot{a}a^{-1}b$ only depends on the 1-forms $B=(\dot{a}a^{-1}, \dot{b}-\dot{a}a^{-1}b)$. We notice that $a^{-1}\dot{a}$ and $\dot{a}a^{-1}$ are both $3\times 3$ skewsymmetric time depending matrices that may be quite different.\\

\noindent
{\bf REMARK  1.6}: A computation in local coordinates for the case of the movement of a rigid body shows that the action of the $3\times 3$ skewsymmetric matrix $\dot{a}a^{-1}$ on the position $x$ at time $t$ just amounts to the vector product by the {\it vortex vector} $\omega=\frac{1}{2}curl (v)$ (See [1,2,3,21] for more details).   \\

The above particular case, well known by anybody studying the analytical mechanics of rigid bodies, can be generalized as follows. If $X$ is a manifold and $G$ is a lie group ({\it not acting necessarily on} $X$), let us consider maps $a:X\rightarrow G: (x)\rightarrow (a(x))$ or equivalently sections of the trivial (principal) bundle $X\times G$ over $X$. If $x+dx$ is a point of $X$ close to $x$, then $T(a)$ will provide a point $a+da=a+\frac{\partial a}{\partial x}dx$ close to $a$ on $G$. We may bring $a$ back to $e$ on $G$ by acting on $a$ with $a^{-1}$, {\it  either on the left or on the right}, getting therefore a $1$-form $a^{-1}da=A$ or $daa^{-1}=B$ with value in $\cal{G}$. As $aa^{-1}=e$ we also get $daa^{-1}=-ada^{-1}=-b^{-1}db$ if we set $b=a^{-1}$ as a way to link $A$ with $B$. When there is an action $y=ax$, we have $x=a^{-1}y=by$ and thus $dy=dax=daa^{-1}y$, a result leading through the first fundamental theorem of Lie to the equivalent formulas:\\
\[   a^{-1}da=A=({A}^{\tau}_i(x)dx^i=-{\omega}^{\tau}_{\sigma}(b(x)){\partial}_ib^{\sigma}(x)dx^i)  \]
\[   daa^{-1}=B=({B}^{\tau}_i(x)dx^i={\omega}^{\tau}_{\sigma}(a(x)){\partial}_ia^{\sigma}(x)dx^i)  \]
Introducing the induced bracket $[A,A](\xi,\eta)=[A(\xi),A(\eta)]\in {\cal{G}}, \forall \xi,\eta\in T$, we may define the $2$-form $dA-[A,A]=F\in {\wedge}^2T^*\otimes {\cal{G}}$ by the local formula ({\it care again to the sign}):\\
\[     {\partial}_iA^{\tau}_j(x)-{\partial}_jA^{\tau}_i(x)-c^{\tau}_{\rho\sigma}A^{\rho}_i(x)A^{\sigma}_j(x)=F^{\tau}_{ij}(x)  \]
This definition can also be adapted to $B$ by using $dB+[B,B]$ and we obtain from the second fundamental theorem of Lie:\\

\noindent
{\bf THEOREM 1.7}: There is a {\it nonlinear gauge sequence}:\\
\[  \begin{array}{ccccc}
X\times G & \longrightarrow & T^*\otimes {\cal{G}} &\stackrel{MC}{ \longrightarrow} & {\wedge}^2T^*\otimes {\cal{G}}  \\
a                & \longrightarrow  &    a^{-1}da=A         &    \longrightarrow & dA-[A,A]=F
\end{array}   \]

Choosing $a$ "close" to $e$, that is $a(x)=e+t\lambda(x)+...$ and linearizing as usual, we obtain the linear operator $d:{\wedge}^0T^*\otimes {\cal{G}}\rightarrow {\wedge}^1T^*\otimes {\cal{G}}:({\lambda}^{\tau}(x))\rightarrow ({\partial}_i{\lambda}^{\tau}(x))$ leading to:\\

\noindent
{\bf COROLLARY  1.8}: There is a {\it linear gauge sequence}:\\ 
\[  {\wedge}^0T^*\otimes {\cal{G}}\stackrel{d}{\longrightarrow} {\wedge}^1T^*\otimes {\cal{G}} \stackrel{d}{\longrightarrow} {\wedge}^2T^*\otimes{\cal{G}} \stackrel{d}{\longrightarrow} ... \stackrel{d}{\longrightarrow} {\wedge}^nT^*\otimes {\cal{G}}\longrightarrow  0   \]
which is the tensor product by $\cal{G}$ of the Poincar\'{e} sequence:\\

 It just remains to introduce the previous results into a variational framework. For this, we may consider a lagrangian on $T^*\otimes \cal{G}$, that is an {\it action} $W=\int w(A)dx$ where $dx=dx^1\wedge ...\wedge dx^n$ and to vary it. With $A=a^{-1}da=-dbb^{-1}$ we may introduce $\lambda=a^{-1}\delta a=-\delta bb^{-1}\in {\cal{G}}={\wedge}^0T^*\otimes {\cal{G}}$ with local coordinates ${\lambda}^{\tau}(x)=-{\omega}^{\tau}_{\sigma}(b(x))\delta b^{\sigma}(x)$ and we obtain $ \delta A=d\lambda - [A,\lambda]  $ that is $\delta A^{\tau}_i={\partial}_i\lambda^{\tau}-c^{\tau}_{\rho\sigma}A^{\rho}_i{\lambda}^{\sigma}$ in local coordinates. Then, setting $\partial w/\partial A={\cal{A}}=({\cal{A}}^i_{\tau})\in {\wedge}^{n-1}T^*\otimes \cal{G}$, we get:\\
\[  \delta W=\int {\cal{A}}\delta Adx=\int {\cal{A}}(d\lambda-[A,\lambda])dx  \]
and therefore, after integration by part, the Euler-Lagrange (EL) equations [13,16,17]:\\
\[     {\partial}_i{\cal{A}}^i_{\tau}+c^{\sigma}_{\rho\tau}A^{\rho}_i{\cal{A}}^i_{\sigma}=0    \]
Such a linear operator for $\cal{A}$ has {\it non-constant coefficients} linearly depending on $A$. However, setting $\delta aa^{-1}=\mu\in {\cal{G}}$, we get $\lambda=a^{-1}(\delta aa^{-1})a=Ad(a)\mu$ while, setting $a'=ab$, we get the {\it gauge transformation} $A \rightarrow A'=(ab)^{-1}d(ab)=b^{-1}a^{-1}(dab+adb)=Ad(b)A+b^{-1}db, \forall b\in G$. Setting $b=e+t\lambda+...$ with $t\ll 1$, then $\delta A$ becomes an infinitesimal gauge transformation. Finally, $ a'=ba\Rightarrow A'=a^{-1}b^{-1}(dba+adb)=a^{-1}(b^{-1}db)a+A\Rightarrow \delta A=Ad(a)d\mu$ when $b=e+t\mu +...$ with $t\ll 1$.
Therefore, introducing $\cal{B}$ such that ${\cal{B}}\mu= {\cal{A}}\lambda$, we get the {\it divergence-like equations} $ {\partial}_i{\cal{B}}^i_{\sigma}=0$.  \\

In 1954, at the birth of GT, the above notations were coming from electromagnetism (EM) with EM {\it potential} $A\in T^*$ and EM {\it field} $dA=F\in {\wedge}^2T^*$ in the relativistic Maxwell theory [12]. Accordingly, $G=U(1)$ (unit circle in the complex plane)$\longrightarrow dim ({\cal{G}})=1$ was the {\it only possibility} to get {\it pure} $1$-form $A$ and $2$-form $F$ when $c=0$. However, "{\it surprisingly}", this result is {\it not coherent at all} with elasticity theory and, {\it a fortiori} with the analytical mechanics of rigid bodies where the Lagrangian is a quadratic expression of 1-forms as we saw because the EM lagrangian $(\epsilon/2)E^2-(1/2{\mu})B^2$  is a quadratic expression of the EM field $F$ as a $2$-form satisfying the first set of Maxwell equations $dF=0$. The dielectric constant $\epsilon$ and the magnetic constant $\mu$ are leading to the electric induction $\vec{D}=\epsilon\vec{E}$ and the magnetic induction $\vec{H}=(1/\mu)\vec{B}$ in the second set of Maxwell equations. In view of the existence of well known field-matter couplings (piezoelectricity, 
photoelasticity) [16,20,23], {\it such a situation is contradictory} as it should lead to put on equal footing $1$-forms and $2$-forms but no other substitute could have been provided at that time, despite the tentatives of the brothers Eugene Cosserat (1866-1931) and Francois Cosserat (1852-1914) in 1909 [5,16,22,23] or of Herman Weyl (1885-1955) in 1918 [16,29] .\\

After this long introduction, the purpose of this paper will be to escape from such a contradiction by using new mathematical tools coming from the formal theory of systems of PD equations and Lie pseudogroups, exactly as we did in [24] for general relativity (GR). In particular, the titles of the three parts that follow will be quite similar to those of this reference though, of course, the contents will be different. The first part proves hat the name "{\it curvature}" given to $F$ has been quite misleading, the resulting confusion between translation and rotation being presented with humour in [32] through the chinese saying " {\it to put Chang's cap on Li's head} ". The second part explains why the Cosserat/Maxwell/Weyl (CMW) theory MUST be described by the Spencer sequence and NOT by the Janet sequence, with a SHIFT by one step contradicting the mathematical foundations of both GR and GT. The third part finally presents the Poincare duality scheme of physics  by means of unexpected methods of homological algebra and algebraic analysis.   \\
 
\noindent
{\bf 2. First Part: The Nonlinear Janet and Spencer Sequences} \\

In 1890, Lie discovered that {\it Lie groups of transformations} were examples of {\it Lie pseudogroups of transformations} along the following definition:  \\
 
 \noindent
{\bf DEFINITION 2.1}: A {\it Lie pseudogroup of transformations} $\Gamma\subset aut(X)$ is a group of transformations solutions of a system of OD or PD equations such that, if $y=f(x)$ and $z=g(y)$ are two solutions, called {\it finite transformations}, that can be composed, then $z=g\circ f(x)=h(x)$ and $x=f^{-1}(y)=g(y)$ are also solutions while $y=x$ is the {\it identity} solution denoted by $id=id_X$ and we shall set $id_q=j_q(id)$.  In all the sequel we shall suppose that $\Gamma$ is {\it transitive} that is $\forall x,y\in X, \exists f\in \Gamma, y=f(x)$ \\

From now on, we shall use the same notations and definitions as in [17,23,24] for jet bundles. In particular, we recall that, if $J_q({\cal{E}})\rightarrow X:(x,y_q)\rightarrow (x) $ is the $q$-jet bundle of ${\cal{E}}\rightarrow X:(x,y)\rightarrow (x) $ with local coordinates $(x^i,y^k_{\mu})$ for $i=1,...,n$, $k=1,...,m$, $0\leq \mid \mu \mid\leq q$ and $y^k_0=y^k$, we may consider sections $f_q:(x)\rightarrow (x,f^k(x), f^k_i(x), f^k_{ij}(x), ...)=(x,f_q(x))$ transforming like the sections $j_q(f):(x) \rightarrow (x,f^k(x), {\partial}_if^k(x), {\partial}_{ij}f^k(x), ...)=(x,j_q(f)(x))$ where both $f_q$ and $j_q(f)$ are over the section $f:(x)\rightarrow (x,y^k=f^k(x))=(x,f(x))$ of $\cal{E}$. The (nonlinear) {\it Spencer operator} just allows to distinguish a section $f_q$ from a section $j_q(f)$ by introducing a kind of "{\it difference}" through the operator $D:J_{q+1}({\cal{E}})\rightarrow T^*\otimes V(J_q({\cal{E}})): f_{q+1}\rightarrow j_1(f_q)-f_{q+1}$ with local components $({\partial}_if^k(x)-f^k_i(x), {\partial}_if^k_j(x)-f^k_{ij}(x),...) $ and more generally $(Df_{q+1})^k_{\mu,i}(x)={\partial}_if^k_{\mu}(x)-f^k_{\mu+1_i}(x)$. If $m=n$ and ${\cal{E}}=X\times X$ with source projection, we denote by ${\Pi}_q={\Pi}_q(X,X)\subset J_q(X\times X)$ the open sub-bundle locally defined by $det(y^k_i)\neq 0$.  \\
 
We also notice that an action $y=f(x,a)$ provides a Lie pseudogroup by eliminating the $p$ parameters $a$ among the equations $y_q=j_q(f)(x,a)$ obtained by successive differentiations with respect to $x$ only when $q$ is large enough. The system ${\cal{R}}_q\subset {\Pi}_q$ of PD equations thus obtained may be quite nonlinear and of high order. Looking for transformations "close" to the identity, that is setting $y=x+t\xi(x)+...$ when $t\ll 1$ is a small constant parameter and passing to the limit $t\rightarrow 0$, we may linearize the above (nonlinear) {\it system of finite Lie equations} in order to obtain a (linear) {\it system of infinitesimal Lie equations} $R_q=id^{-1}_q(V({\cal{R}}_q))\subset J_q(T)$ for vector fields. Such a system has the property that, if $\xi,\eta$ are two solutions, then $[\xi,\eta]$ is also a solution. Accordingly, the set $\Theta\subset T$ of its solutions satisfies $[\Theta,\Theta]\subset \Theta$ and can therefore be considered as the Lie algebra of $\Gamma$.\\

\noindent
{\bf GAUGING PROCEDURE REVISITED 2.2 } : When there is a Lie group of transformations, setting $f(x)=f(x,a(x))$ and $f_q(x)=j_q(f)(x,a(x))$, we obtain $a(x)=a=cst \Leftrightarrow f_q=j_q(f)$ because $Df_{q+1}=j_1(f_q)-f_{q+1}=(\partial f_q(x,a(x))/\partial a^{\tau}){\partial}_ia^{\tau}(x)$ and the matrix involved has rank $p$ in the following commutative diagram:  \\
\[  \begin{array}{rcccl}
0\rightarrow  & \hspace{5mm} X  \times  G &  =  &  \hspace{4mm}{\cal{R}}_q  &  \rightarrow 0  \\
      &   a=cst \uparrow \downarrow \uparrow a(x)  &    &  j_q(f) \uparrow \downarrow \uparrow  f_q &  \\
         &  \hspace{4mm}X     &  =   &    \hspace{4mm}X   &   
\end{array}  \]

Looking at the way a vector field and its derivatives are transformed under any $f\in aut(X)$ while replacing $j_q(f)$ by $f_q$, we obtain:\\
\[  {\eta}^k(f(x))=f^k_r(x){\xi}^r(x) \Rightarrow {\eta}^k_u(f(x))f^u_i(x)=f^k_r(x){\xi}^r_i(x)+f^k_{ri}(x){\xi}^r(x)\]
and so on, a result leading to:\\

\noindent
{\bf LEMMA 2.3}: $J_q(T)$ is {\it associated} with ${\Pi}_{q+1}$ that is we can obtain a new section ${\eta}_q=f_{q+1}({\xi}_q)$ from any section ${\xi}_q \in J_q(T)$ and any section $f_{q+1}\in {\Pi}_{q+1}$ by the formula:\\
\[ d_{\mu}{\eta}^k\equiv {\eta}^k_rf^r_{\mu}+ ...=f^k_r{\xi}^r_{\mu}+  ...  +f^k_{\mu+1_r}{\xi}^r , \forall 0\leq {\mid}\mu {\mid}\leq q\]
where the left member belongs to $V({\Pi}_q)$. Similarly $R_q\subset J_q(T)$ is associated with ${\cal{R}}_{q+1}\subset {\Pi}_{q+1}$.\\

In order to construct another nonlinear sequence, we need a few basic definitions on {\it Lie groupoids} and {\it Lie algebroids} that will become substitutes for Lie groups and Lie algebras. The first idea is to use the chain rule for derivatives $j_q(g\circ f)=j_q(g)\circ j_q(f)$ whenever $f,g\in aut(X)$ can be composed and to replace both $j_q(f)$ and $j_q(g)$ respectively by $f_q$ and $g_q$ in order to obtain the new section $g_q\circ f_q$. This kind of "{\it composition}" law can be written in a pointwise symbolic way by introducing another copy $Z$ of $X$ with local coordinates $(z)$ as follows:\\
\[ {\gamma}_q:{\Pi}_q(Y,Z){\times}_Y{\Pi}_q(X,Y)\rightarrow {\Pi}_q(X,Z):((y,z,\frac{\partial z}{\partial y},...),(x,y,\frac{\partial y}{\partial x},...)\rightarrow (x,z,\frac{\partial z}{\partial y}\frac{\partial y}{\partial x},...)      \]
We may also define $j_q(f)^{-1}=j_q(f^{-1})$ and obtain similarly an "{\it inversion}" law.\\

\noindent
{\bf DEFINITION 2.4}: A fibered submanifold ${\cal{R}}_q\subset {\Pi}_q$ is called a {\it system of finite Lie equations} or a {\it Lie groupoid} of order $q$ if we have an induced {\it source projection} ${\alpha}_q:{\cal{R}}_q\rightarrow X$, {\it target projection} ${\beta}_q:{\cal{R}}_q\rightarrow X$, {\it composition} ${\gamma}_q:{\cal{R}}_q{\times}_X{\cal{R}}_q\rightarrow {\cal{R}}_q$, {\it inversion} ${\iota}_q:{\cal{R}}_q\rightarrow {\cal{R}}_q$ 
and {\it identity} $id_q:X\rightarrow {\cal{R}}_q$. In the sequel we shall only consider {\it transitive} Lie groupoids such that the map $({\alpha}_q,{\beta}_q):{\cal{R}}_q\rightarrow X\times X $ is an epimorphism. One can prove that the new system ${\rho}_r({\cal{R}}_q)={\cal{R}}_{q+r}$ obtained by differentiating $r$ times all the defining equations of ${\cal{R}}_q$ is a Lie groupoid of order $q+r$. \\

Now, using the {\it algebraic bracket} $\{ j_{q+1}(\xi),j_{q+1}(\eta)\}=j_q([\xi,\eta]), \forall \xi,\eta\in T$, we may  obtain by bilinearity a {\it differential bracket} on $J_q(T)$ extending the bracket on $T$:\\
\[   [{\xi}_q,{\eta}_q]=\{{\xi}_{q+1},{\eta}_{q+1}\}+i(\xi)D{\eta}_{q+1}-i(\eta)D{\xi}_{q+1}, \forall {\xi}_q,{\eta}_q\in J_q(T) \]
which does not depend on the respective lifts ${\xi}_{q+1}$ and ${\eta}_{q+1}$ of ${\xi}_q$ and ${\eta}_q$ in $J_{q+1}(T)$. One can prove that his bracket on sections satisfies the Jacobi identity and we set: \\

\noindent
{\bf DEFINITION 2.5}: We say that a vector subbundle $R_q\subset J_q(T)$ is a {\it system of infinitesimal Lie equations} or a {\it Lie algebroid} if $[R_q,R_q]\subset R_q$, that is to say $[{\xi}_q,{\eta}_q]\in R_q, \forall {\xi}_q,{\eta}_q\in R_q$. Such a definition can be tested by means of computer algebra.  \\

\noindent
{\bf EXAMPLE 2.6}: With $n=1, q=2, X=\mathbb{R}$ and evident notations, the components of $[{\xi}_2,{\eta}_2]$ at order zero, one and two are defined by the totally unusual successive formulas:\\
\[    [\xi,\eta]=\xi{\partial}_x\eta-\eta{\partial}_x\xi     \]
\[    ([{\xi}_1,{\eta}_1])_x=\xi{\partial}_x{\eta}_x-\eta{\partial}_x{\xi}_x    \]
\[    ([{\xi}_2,{\eta}_2])_{xx}={\xi}_x{\eta}_{xx}-{\eta}_x{\xi}_{xx}+\xi{\partial}_x{\eta}_{xx}-\eta{\partial}_x{\xi}_{xx}   \]
For affine transformations, ${\xi}_{xx}=0,{\eta}_{xx}=0\Rightarrow ([{\xi}_2,{\eta}_2])_{xx}=0$ and thus $[R_2,R_2]\subset R_2$.\\

We may {\it prolong} the vertical infinitesimal transformations $\eta={\eta}^k(y)\frac{\partial}{\partial y^k}$ to the jet coordinates up to order $q$ in order to obtain:\\
\[   {\eta}^k(y)\frac{\partial}{\partial y^k}+\frac{\partial {\eta}^k}{\partial y^r}y^r_i\frac{\partial}{\partial y^k_i}+(\frac{{\partial}^2{\eta}^k}{\partial y^r\partial y^s}y^r_iy^s_j+\frac{\partial {\eta}^k}{\partial y^r}y^r_{ij})\frac{\partial}{\partial y^k_{ij}}+...    \]
where we have replaced $j_q(f)(x)$ by $y_q$, each component beeing the "formal" derivative of the previous one. Replacing $j_q(\eta)$ by ${\eta}_q$ as sections of $R_q$ over the target, we obtain a vertical vector field $\sharp ({\eta}_q)$ over ${\Pi}_q$ such that $[\sharp({\eta}_q),\sharp ({\zeta}_q)]=\sharp ([{\eta}_q,{\zeta}_q]), \forall {\eta}_q,{\zeta}_q\in R_q$ over the target. We may then use the Frobenius theorem in order to find a generating fundamental set of {\it differential invariants} $\{{\Phi}^{\tau}(y_q)\}$ up to order $q$ which are such that ${\Phi}^{\tau}({\bar{y}}_q)={\Phi}^{\tau}(y_q)$ by using the chain rule for derivatives whenever $\bar{y}=g(y)\in \Gamma$ acting now on $Y$. Looking at the way the differential invariants are transformed between themselves under changes of source, we may define a natural bundle ${\cal{F}}\rightarrow X:(x,u)\rightarrow (x)$. Specializing the ${\Phi}^{\tau}$ at $id_q(x)$ we obtain the {\it Lie form} ${\Phi}^{\tau}(y_q)={\omega}^{\tau}(x)$ of ${\cal{R}}_q$ and a section $\omega:(x)\rightarrow (x,\omega (x))$ of ${\cal{F}}$. If we introduce the maximum number of formal derivatives $d_i{\Phi}^{\tau}$ that are linearly independent over the jets of strict order $q+1$, any other formal derivative is a linear combination with coefficients functions of $y_q$. Applying $\sharp (R_q)$, we get a contradiction unless these coefficients are killed by $\sharp (R_q)$ and are thus functions of the fundamental set, a result leading to CC of the form $ I(j_1(\omega))\equiv A(\omega){\partial}_x\omega+B(\omega)=0$. Finally, setting $v=A(u)u_x+B(u)$, we obtain a new natural bundle ${\cal{F}}_1\rightarrow X:(x,u,v)\rightarrow ((x)$ as a vector bundle over ${\cal{F}}$.   \\

 \noindent
{\bf THEOREM 2.7}: There exists a {\it nonlinear Janet sequence} associated with the Lie form of an involutive system of finite Lie equations:   \\
\[  \begin{array}{rcccl}
  & \Phi \circ j_q &   &  I \circ j_1  &   \\
  0\rightarrow \Gamma \rightarrow aut(X) &\rightrightarrows &  {\cal{F}}  &\rightrightarrows   &  {\cal{F}}_1 \\
     & \omega\circ\alpha  &  &  0 &
  \end{array}  \]
where the kernel of the first operator $f\rightarrow \Phi \circ j_q(f)= \Phi (j_q(f))=j_q(f)^{-1}(\omega)$ is taken with respect to the section $\omega$ of $\cal{F}$ while the kernel of the second operator is taken with respect to the zero section of the vector bundle ${\cal{F}}_1$ over ${\cal{F}}$ (Compare to [10,28]).\\

 \noindent
{\bf THEOREM 2.8} : There is a {\it first nonlinear Spencer sequence}:\\
\[ 0\longrightarrow aut(X) \stackrel{j_{q+1}}{\longrightarrow} {\Pi}_{q+1}(X,X)\stackrel{\bar{D}}{\longrightarrow}T^*\otimes J_q(T)\stackrel{{\bar{D}}'}{\longrightarrow} {\wedge}^2T^*\otimes J_{q-1}(T)  \]
with $\bar{D}f_{q+1}\equiv f_{q+1}^{-1}\circ j_1(f_q)-id_{q+1}={\chi}_q \Rightarrow {\bar{D}}'{\chi}_q(\xi,\eta)\equiv D{\chi}_q(\xi,\eta)-\{{\chi}_q(\xi),{\chi}_q(\eta)\}=0 $. Moreover, setting ${\chi}_0=A-id\in T^*\otimes T$, this sequence is locally exact if $det(A)\neq 0$ and there is an induced  {\it second nonlinear Spencer sequence} (See next section for the definition of the Spencer bundles): \\
\[  0 \longrightarrow aut(X) \stackrel{j_q}{\longrightarrow} {\cal{R}}_q \stackrel{{\bar{D}}_1}{\longrightarrow} C_1(T) \stackrel{{\bar{D}}_2}{\longrightarrow} C_2(T)  \]
where all the operators involved are involutive and $C_1(T), C_2(T)$ linearly depend on $J_q(T)$ {\it only}.  \\

\noindent
{\bf Proof}: There is a canonical inclusion ${\Pi}_{q+1}\subset J_1({\Pi}_q)$ defined by $y^k_{\mu,i}=y^k_{\mu+1_i}$ and the composition $f^{-1}_{q+1}\circ j_1(f_q)$ is a well defined section of $J_1({\Pi}_q)$ over the section $f^{-1}_q\circ f_q=id_q$ of ${\Pi}_q$ like $id_{q+1}$. The difference ${\chi}_q=f^{-1}_{q+1}\circ j_1(f_q)-id_{q+1}$ is thus a section of $T^*\otimes V({\Pi}_q)$ over $id_q$ and we have already noticed that 
$ id^{-1}_q(V({\Pi}_q))=J_q(T)$. For $q=1$ we get with $g_1=f^{-1}_1$:\\
\[ {\chi}^k_{,i}=g^k_l{\partial}_if^l-{\delta}^k_i=A^k_i-{\delta}^k_i,\hspace{5mm} {\chi}^k_{j,i}=g^k_l({\partial}_if^l_j-A^r_if^l_{rj})  \]
We also obtain from Lemma 2.3 the useful formula $ f^k_r{\chi}^r_{\mu,i}+...+f^k_{\mu+1_r}{\chi}^r_{,i}={\partial}_if^k_{\mu}-f^k_{\mu+1_i}$ allowing to determine ${\chi}_q$ inductively with ${\chi}^k_{\mu,i}= -g^k_lA^r_if^l_{\mu+1_r} + (order \leq \mid\mu\mid)$ when $q \geq 1$. It just remains to set ${\chi}_q={\tau}_q\circ A$ as $1$-forms in order to construct $C_1$ and $C_2$ by quotients.  \\
We refer to ([17], p 215) for the inductive proof of the local exactness, providing the only formulas that will be used later on and can be checked directly by the reader:\\
\[  {\partial}_i{\chi}^k_{,j}-{\partial}_j{\chi}^k_{,i}-{\chi}^k_{i,j}+{\chi}^k_{j,i}-({\chi}^r_{,i}{\chi}^k_{r,j}-{\chi}^r_{,j}{\chi}^k_{r,i})=0  \]
\[  {\partial}_i{\chi}^k_{l,j}-{\partial}_j{\chi}^k_{l,i}-{\chi}^k_{li,j}+{\chi}^k_{lj,i}-({\chi}^r_{,i}{\chi}^k_{lr,j}+{\chi}^r_{l,i}{\chi}^k_{r,j}-{\chi}^r_{l,j}{\chi}^k_{r,i}-{\chi}^r_{,j}{\chi}^k_{lr,i})=0\]
There is no need for double-arrows in this framework as the kernels are taken with respect to the zero section of the vector bundles involved. We finally notice that the main difference with the gauge sequence is that {\it all the indices range from} $1$ {\it to} $n$ and that the condition $det(A)\neq 0$ amounts to $\Delta=det({\partial}_if^k)\neq 0$ because $det(f^k_i)\neq 0$ by assumption (See [15,17,23] for more details).Ê \\
\hspace*{12cm}  Q.E.D.  \\

\noindent
{\bf COROLLARY 2.9}: There is a {\it first restricted nonlinear Spencer sequence}:\\
\[ 0\longrightarrow \Gamma \stackrel{j_{q+1}}{\longrightarrow} {\cal{R}}_{q+1} \stackrel{\bar{D}}{\longrightarrow} T^*\otimes R_q\stackrel{{\bar{D}}'}{\longrightarrow}{\wedge}^2T^*\otimes J_{q-1}(T)  \]
and an induced {\it second restricted nonlinear Spencer sequence}:  \\
\[ 0 \longrightarrow \Gamma \stackrel{j_q}{\longrightarrow} {\cal{R}}_q \stackrel{{\bar{D}}_1}{\longrightarrow} C_1 \stackrel{{\bar{D}}_2}{\longrightarrow} C_2 \]
where all the operators involved are involutive and $C_1,C_2$ linearly depend on $R_q$ {\it only}. This sequence is locally isomorphic to the corresponding gauge sequence for any Lie group of transformations when $q$ is large enough. {\it The action, which is essential in the Spencer sequence, disappears in the gauge sequence}.   \\

\noindent
{\bf DEFINITION 2.10}: A {\it splitting} of the short exact sequence $0\rightarrow R^0_q\rightarrow R_q\stackrel{{\pi}^q_0}{\rightarrow} T \rightarrow 0$ is a map ${\chi}'_q:T\rightarrow R_q$ such that ${\pi}^q_0\circ {\chi}'_q=id_T$ or equivalently a section of $T^*\otimes R_q$ over $id_T\in T^*\otimes T$ and is called a $R_q$-{\it connection}. Its {\it curvature} ${\kappa}'_q\in {\wedge}^2T^*\otimes R^0_q$ is defined by ${\kappa}'_q(\xi,\eta)=[{\chi}'_q(\xi),{\chi}'_q(\eta)]-{\chi}'_q([\xi,\eta])$. We notice that ${\chi}'_q=-{\chi}_q$ is a connection with ${\bar{D}}'{\chi}'_q={\kappa}'_q$ {if and only if} $A=0$ but {\it connections cannot be used for describing fields because we must have} $\Delta \neq 0$. {\. In particular $({\delta}^k_i,-{\gamma}^k_{ij})$ is the only existing symmetric connection for the Killing system.         \\

\noindent
{\bf REMARK 2.11}: Rewriting the previous formulas with $A$ instead of ${\chi}_0$ we get:  \\
\[ {\partial}_iA^k_j-{\partial}_jA^k_i-A^r_i{\chi}^k_{r,j}+A^r_j{\chi}^k_{r,i}=0  \]
\[ {\partial}_i{\chi}^k_{l,j}-{\partial}_j{\chi}^k_{l,i}-{\chi}^r_{l,i}{\chi}^k_{r,j}+{\chi}^r_{l,j}{\chi}^k_{r,i}-A^r_i{\chi}^k_{lr,j}+A^r_j{\chi}^k_{lr,i}=0  \]
When $q=1, g_2=0$ and though surprising it may look like, we find back {\it exactly} all the formulas presented by E. and F. Cosserat in ([5], p 123 and [27]) (Compare to [10]). \\
 
Finally, setting $f'_{q+1}=g_{q+1}\circ f_{q+1}$, we get $\bar{D}f'_{q+1}=f^{-1}_{q+1}\circ g^{-1}_{q+1}\circ j_1(g_q)\circ j_1(f_q)-id_{q+1}=f^{-1}_{q+1}\circ \bar{D}g_{q+1}\circ j_1(f_q)+\bar{D}f_{q+1}, \forall f_{q+1},g_{q+1}\in {\cal{R}}_{q+1}$. With ${\chi}_q=\bar{D}g_{q+1}$, we get the {\it gauge transformation} ${\chi}_q \rightarrow f^{-1}_{q+1}\circ {\chi}_q\circ j_1(f_q)+\bar{D}f_{q+1}, \forall f_{q+1}\in {\cal{R}}_{q+1}$ as in the introduction, thus ACTING ON THE FIELDS ${\chi}_q$ WHILE PRESERVING THE FIELD EQUATIONS ${\bar{D}}'{\chi}_q=0$. Setting $f_{q+1}=id_{q+1}+t{\xi}_{q+1}+...$ with $t\ll 1$ {\it over the source}, we obtain an {\it infinitesimal gauge transformation} of the form $\delta {\chi}_q=D{\xi}_{q+1} +L(j_1({\xi}_{q+1}){\chi}_q$ as in [16,17,23]. However, setting now ${\chi}_q=\bar{D}f_{q+1}$ and $g_{q+1}=id_{q+1}+t{\eta}_{q+1}+...$ with $t\ll 1$ {\it over the target}, we get $\delta {\chi}_q=f^{-1}_{q+1}\circ D{\eta}_{q+1}\circ j_1(f_q)$. The same variation is obtained whenever ${\eta}_{q+1}=f_{q+2}({\xi}_{q+1}+{\chi}_{q+1}(\xi))$ with ${\chi}_{q+1}=\bar{D}f_{q+2}$, a transformation which only depends on $j_1(f_{q+1})$ and is invertible if and only if $det(A)\neq 0$ [16,17]. This result proves that $J_q(T)$ is also associated with the groupoid ${\Pi}_{q,1} \subset J_1({\Pi}_q)$ defined by $det(y^k_{0,i}) \neq 0$. With $g_1=f^{-1}_1$, we have the unusual formulas: \\
\[   {\eta}^k={\xi}^r{\partial}_rf^k, \hspace{3mm} {\eta}^k_u=g^i_uf^k_r{\xi}^r_i+g^i_u{\xi}^r{\partial}_rf^k_i . \] 
Accordingly, THE DUAL EQUATIONS WILL ONLY DEPEND ON THE LINEAR SPENCER OPERATOR $D$. Moreover, in view of the two variational results obtained at the end of the introduction, THE CMW EQUATIONS CANNOT COME FROM THE GAUGE SEQUENCE, contrary to what mechanicians still believe after more than a century.  \\ 

\noindent
{\bf EXAMPLE 2.12}: We have the formulas (Compare to [5] and [29],(76) p 289,(78) p 290):  \\
\[ \delta {\chi}^k_{,i}=({\partial}_i{\xi}^k-{\xi}^k_i)+({\xi}^r{\partial}_r{\chi}^k_{,i}+{\chi}^k_{,r}{\partial}_i{\xi}^r-{\chi}^r_{,i}{\xi}^k_r)=g^k_v(\frac{\partial {\eta}^v}{\partial y^u} -{\eta}^v_u){\partial}_if^u   \]
\[  \delta {\chi}^k_{j,i}=({\partial}_i{\xi}^k_j-{\xi}^k_{ij})+({\xi}^r{\partial}_r{\chi}^k_{j,i}+{\chi}^k_{j,r}{\partial}_i{\xi}^r+{\chi}^k_{r,i}{\xi}^r_j-{\chi}^r_{j,i}{\xi}^k_r-{\chi}^r_{,i}{\xi}^k_{jr} )  \]
Setting ${\alpha}_i={\chi}^r_{r,i}$, we have $\delta{\alpha}_i={(\partial}_i{\xi}^r_r-{\xi}^r_{ri})+({\xi}^r{\partial}_r{\alpha}_i+{\alpha}_r{\partial}_i{\xi}^r-{\chi}^s_{,i}{\xi}^r_{rs})  $.  \\

\noindent
{\bf EXAMPLE 2.13}: (Projective transformations) With ${\xi}_{xxx}=0$, the formal adjoint of the Spencer operator brings as many dual equations as the number of parameters (1 translation + 1 dilatation + 1 elation). \\  
\[ \begin{array}{rcl}
  \sigma ({\partial}_x\xi-{\xi}_x)+\mu({\partial}_x{\xi}_x-{\xi}_{xx})+\nu ({\partial}_x{\xi}_{xx}-{\xi}_{xxx})&=& -[ ({\partial}_x\sigma) \xi + ({\partial}_x\mu +\sigma){\xi}_x + ({\partial}_x\nu + \mu){\xi}_{xx}] \\
  &  &  + {\partial}_x(\sigma \xi +\mu {\xi}_x + \nu {\xi}_{xx})
  \end{array}  \]
\noindent
{\it Cosserat/Weyl equations} :  $ {\partial}_x\sigma=f  \hspace{1mm}, \hspace{2mm} {\partial}_x\mu + \sigma = m \hspace{1mm},\hspace{2mm} {\partial}_x\nu +\mu=j $      \hspace{2mm} (equivalent "momenta") \\
  
\noindent 
{\bf 3. Second Part: The Linear Janet and Spencer Sequences}  \\ 

It remains to understand how the shift by one step in the interpretation of the Spencer sequence is coherent with mechanics and electromagnetism both with their well known couplings [16,17,23,24]. In a word, the problem we have to solve is to get a $2$-form in ${\wedge}^2T^*$ from a $1$-form in $T^*\otimes R_q$ for a certain $R_q\subset J_q(T)$.\\

For this purpose, introducing the {\it Spencer map} $\delta:{\wedge}^sT^*\otimes S_{q+1}T^*\otimes E\rightarrow {\wedge}^{s+1}T^*\otimes S_qT^*\otimes E$ defined by $({\delta \omega})^k_{\mu}=dx^i\wedge {\omega}^k_{\mu+1_i}$, we recall from [17,24,26] the definition of the {\it Janet bundles} $ F_r={\wedge}^rT^*\otimes J_q(E)/({\wedge}^rT^*\otimes R_q + \delta ({\wedge}^{r-1}T^*\otimes S_{q+1}T^*\otimes E)) $ and the {\it Spencer bundles} $C_r={\wedge}^rT^*\otimes R_q/{\delta}({\wedge}^{r-1}T^*\otimes g_{q+1})$ or $C_r(E)={\wedge}^rT^*\otimes J_q(E)/\delta({\wedge}^{r-1}T^*\otimes S_{q+1}T^*\otimes  E)$ with $C_r\subset C_r(E)$. When $R_q\subset J_q(E)$ is an involutive system on $E$, we have the following crucial commutative diagram with exact columns where each operator involved is first order apart from ${\cal{D}}=\Phi\circ j_q$, generates the CC of the preceding one and is induced by the extension $D:{\wedge}^rT^*\otimes J_{q+1}(E) \rightarrow {\wedge}^{r+1}T^*\otimes J_q(E):\alpha\otimes {\xi}_{q+1}\rightarrow d\alpha\otimes {\xi}_q+(-1)^r\alpha\wedge D{\xi}_{q+1}$ of the Spencer operator $D:J_{q+1}(E) \rightarrow T^*\otimes J_q(E):{\xi}_{q+1}\rightarrow j_1({\xi}_q)-{\xi}_{q+1}$. The upper sequence is the (second) {\it linear Spencer sequence} while the lower sequence is the {\it linear Janet sequence} and the sum $dim(C_r)+dim(F_r)=dim(C_r(E))$ does not depend on the system while the epimorphisms ${\Phi}_r$ are induced by $\Phi={\Phi}_0$.  \\
  \\
  \[ \begin{array}{rcccccccccccl}
 &&&&& 0 &&0&&0&  &0&  \\
 &&&&& \downarrow && \downarrow && \downarrow &    & \downarrow &  \\
  & 0& \rightarrow& \Theta &\stackrel{j_q}{\longrightarrow}&C_0 &\stackrel{D_1}{\longrightarrow}& C_1 &\stackrel{D_2}{\longrightarrow} & C_2 &\stackrel{D_3}{\longrightarrow} ... \stackrel{D_n}{\rightarrow}& C_n &\rightarrow 0 \\
  &&&&& \downarrow & & \downarrow & & \downarrow & &\downarrow &     \\
   & 0 & \rightarrow & E & \stackrel{j_q}{\longrightarrow} & C_0(E) & \stackrel{D_1}{\longrightarrow} & C_1(E) &\stackrel{D_2}{\longrightarrow} & C_2(E) &\stackrel{D_3}{\longrightarrow} ... \stackrel{D_n}{\longrightarrow} & C_n(E) &   \rightarrow 0 \\
   & & & \parallel && \hspace{5mm}\downarrow {\Phi}_0 & &\hspace{5mm} \downarrow {\Phi}_1 & & \hspace{5mm}\downarrow {\Phi}_2 &  & \hspace{5mm}\downarrow {\Phi}_n & \\
   0 \rightarrow & \Theta &\rightarrow & E & \stackrel{\cal{D}}{\longrightarrow} & F_0  & \stackrel{{\cal{D}}_1}{\longrightarrow} & F_1 & \stackrel{{\cal{D}}_2}{\longrightarrow} & F_2 & \stackrel{{\cal{D}}_3}{\longrightarrow} ... \stackrel{{\cal{D}}_n}{\longrightarrow} & F_n & \rightarrow  0 \\
   &&&&& \downarrow & & \downarrow & & \downarrow &   &\downarrow &   \\
   &&&&& 0 && 0 && 0 &&0 &  
   \end{array}     \]
 For later computations, the sequence $J_3(E)\stackrel{D}{\longrightarrow}T^*\otimes J_2(E)\stackrel{D}{\longrightarrow} {\wedge}^2T^*\otimes J_1(E)$ can be described by the images ${\partial}_i{\xi}^k-{\xi}^k_i=X^k_{,i}$ , ${\partial}_i{\xi}^k_j-{\xi}^k_{ij}=X^k_{j,i}$ , $ {\partial}_i{\xi}^k_{lj}-{\xi}^k_{lij}=X^k_{lj,i}$ leading to the identities:\\
 \[  {\partial}_iX^k_{,j}-{\partial}_jX^k_{,i}+X^k_{j,i}-X^k_{i,j}=0, \hspace{5mm} {\partial}_iX^k_{l,j}-{\partial}_jX^k_{l,i}+X^k_{lj,i}-X^k_{li,j}=0  \]
We also recall that the linear Spencer sequence for a Lie group of transformations $G\times X\rightarrow X$, which {\it essentially} depends on the action because infinitesimal generators are needed, is locally isomorphic to the linear gauge sequence which {\it does not} depend on the action any longer as it is the tensor product of the Poincar\'{e} sequence by the Lie algebra ${\cal{G}}$.  \\

The main idea will be to {\it introduce and compare} the three Lie groups of transformations:ÊÊ\\
\noindent
$\bullet$ The {\it Poincare group} of transformations with $10$ parameters leading to the {\it Killing system} $R_2$:  \\
\[  (L({\xi}_1)\omega)_{ij}\equiv {\omega}_{rj}(x){\xi}^r_i+{\omega}_{ir}(x){\xi}^r_j+{\xi}^r{\partial}_r{\omega}_{ij}(x)=0 \]
\[  (L({\xi}_2)\gamma)^k_{ij}\equiv {\xi}^k_{ij}+{\gamma}^k_{rj}(x){\xi}^r_i+{\gamma}^k_{ir}(x){\xi}^r_j-{\gamma}^r_{ij}(x){\xi}^k_r+{\xi}^r{\partial}_r{\gamma}^k_{ij}(x)=0 \]
\noindent 
$\bullet$ The {\it Weyl group} of transformations with $11$ parameters leading to the system ${\tilde{R}}_2$:   \\
\[  (L({\xi}_1)\omega)_{ij}\equiv {\omega}_{rj}(x){\xi}^r_i+{\omega}_{ir}(x){\xi}^r_j+{\xi}^r{\partial}_r{\omega}_{ij}(x)=A(x){\omega}_{ij}(x) \]
\[  (L({\xi}_2)\gamma)^k_{ij}\equiv {\xi}^k_{ij}+{\gamma}^k_{rj}(x){\xi}^r_i+{\gamma}^k_{ir}(x){\xi}^r_j-{\gamma}^r_{ij}(x){\xi}^k_r+{\xi}^r{\partial}_r{\gamma}^k_{ij}(x)=0 \]
\noindent
$\bullet$ The {\it conformal group} of transformations with $15$ parameters leading to the {\it conformal Killing system} ${\hat{R}}_2$ and to the corresponding Janet/Spencer diagram: \\
\[ (    L({\xi}_1)\omega)_{ij}\equiv {\omega}_{rj}(x){\xi}^r_i+{\omega}_{ir}(x){\xi}^r_j+{\xi}^r{\partial}_r{\omega}_{ij}(x)=A(x){\omega}_{ij}(x) \]
\[ (L({\xi}_2)\gamma)^k_{ij}\equiv {\xi}^k_{ij}+{\gamma}^k_{rj}(x){\xi}^r_i+{\gamma}^k_{ir}(x){\xi}^r_j-{\gamma}^r_{ij}(x){\xi}^k_r+{\xi}^r{\partial}_r{\gamma}^k_{ij}(x)={\delta}^k_iA_j(x)+{\delta}^k_jA_i(x)-{\omega}_{ij}(x){\omega}^{kr}(x)A_r(x) \]
where one has to eliminate the arbitrary function $A(x)$ and $1$-form $A_i(x)dx^i$ for finding sections, replacing the {\it ordinary Lie derivative} ${\cal{L}}(\xi)$ by the {\it formal Lie derivative} $L({\xi}_q)$, that is replacing $j_q(\xi)$ by ${\xi}_q$ when needed. In these formulas, $\omega \in S_2T^*$ with $det(\omega)\neq 0$ and $j_1(\omega)\simeq (\omega,\gamma)$.\\
   \[ \begin{array}{rccccccccccccccl}
 &&&&& 0 &&0&&0&  &0& & 0 & & \\
 &&&&& \downarrow && \downarrow && \downarrow &    & \downarrow & & \downarrow &  \\
  & 0& \rightarrow& \Theta &\stackrel{j_2}{\rightarrow}&15 &\stackrel{D_1}{\rightarrow}& 60 &\stackrel{D_2}{\rightarrow} & 90 &\stackrel{D_3}{\rightarrow} & 60 &  \stackrel{D_4}{\rightarrow}& 15 &\rightarrow 0   \\
  &&&&& \downarrow & & \downarrow & & \downarrow & &\downarrow &  & \downarrow &   \\
   & 0 & \rightarrow & 4 & \stackrel{j_2}{\rightarrow} & 60 & \stackrel{D_1}{\rightarrow} & 160 &\stackrel{D_2}{\rightarrow} & 180 &\stackrel{D_3}{\rightarrow} & 96 & \stackrel{D_4}{\rightarrow} & 20 &   \rightarrow 0 \\
   & & & \parallel && \hspace{5mm}\downarrow {\Phi}_0 & &\hspace{5mm} \downarrow {\Phi}_1 & & \hspace{5mm}\downarrow {\Phi}_2 &  & \hspace{5mm}\downarrow {\Phi}_3 & & \hspace{5mm} \downarrow {\Phi}_4 &  \\
   0 \rightarrow & \Theta &\rightarrow & 4 & \stackrel{\cal{D}}{\rightarrow} & 45  & \stackrel{{\cal{D}}_1}{\rightarrow} & 100 & \stackrel{{\cal{D}}_2}{\rightarrow} & 90 & \stackrel{{\cal{D}}_3}{\rightarrow} & 36 &\stackrel{{\cal{D}}_4}{\rightarrow} & 5 & \rightarrow  0   \\
   &&&&& \downarrow & & \downarrow & & \downarrow &   &\downarrow &  & \downarrow &   \\
   &&&&& 0 && 0 && 0 &&0 &  & 0 &
   \end{array}     \]
We shall use the inclusions $R_2\subset{\tilde{R}}_2\subset{\hat{R}}_2$ in the tricky proof of the next crucial proposition:  \\

\noindent
{\bf PROPOSITION 3.1}: The Spencer sequence for the conformal Lie pseudogroup projects onto the Poincare sequence {\it with a shift by one step}.  \\

\noindent
{\it Proof}: Using $({\delta}^k_i,-{\gamma}^k_{ij})$ as a $R_1$-connection and the fact that $L({\xi}_2)\gamma\in S_2T^*\otimes T, \forall {\xi}_2\in J_2(T)$ while setting $(A^k_{l,i}=X^k_{l,i}+{\gamma}^k_{ls}X^s_{,i}) \in T^*\otimes T^*\otimes T$ with $(A^r_{r,i}=A_i)\in T^*$ and $(B^k_{lj,i}=X^k_{lj,i}+{\gamma}^k_{sj}X^s_{l,i}+{\gamma}^k_{ls}X^s_{j,i}-{\gamma}^s_{lj}X^k_{s,i}+X^r_{,i}{\partial}_r{\gamma}^k_{lj}) \in T^*\otimes S_2T^*\otimes T$ that can be composed with $\delta$ for obtaining the trace, we obtain the following commutative and exact diagram:  \\
\[   \begin{array}{rcccccl}
  &  &  &  0 & & 0 &  \\
  &  &  &  \downarrow & & \downarrow &  \\
    & 0 & \rightarrow & {\hat{g}}_2 & \rightarrow & T^* & \rightarrow 0  \\
    & \downarrow &  & \downarrow & & \parallel &  \\
    0\rightarrow & {\tilde{R}}_2  &  \rightarrow & {\hat{R}}_2 &\rightarrow & T^*  & \rightarrow 0  \\
  & \downarrow &  & \downarrow & & \downarrow &   \\
  0 \rightarrow & {\tilde{R}}_1 & = &{\hat{R}}_1 & \rightarrow  & 0  &  \\
    & \downarrow && \downarrow & & &  \\
    &   0  & &  0 & & & 
    \end{array}     \]
We also obtain from the relations ${\partial}_i{\gamma}^r_{rj}={\partial}_j{\gamma}^r_{ri}$ and the two previous identities:  \\
\[ \begin{array}{rcl}
F_{ij}=B^r_{ri,j}-B^r_{rj,i} & = & X^r_{ri,j}-X^r_{rj,i}+{\gamma}^r_{rs}X^s_{i,j}-{\gamma}^r_{rs}X^s_{j,i}+X^r_{,j}{\partial}_r{\gamma}^s_{si}-X^r_{,i}{\partial}_r{\gamma}^s_{sj}  \\
  &  =  & {\partial}_iX^r_{r,j}-{\partial}_jX^r_{r,i}+{\gamma}^r_{rs}(X^s_{i,j}-X^s_{j,i})+X^r_{,j}{\partial}_i{\gamma}^s_{sr}-X^r_{,i}{\partial}_j{\gamma}^s_{sr} \\
    &  =  & {\partial}_i(X^r_{r,j}+{\gamma}^r_{rs}X^s_{,j})-{\partial}_j(X^r_{r,i}+{\gamma}^r_{rs}X^s_{s,i})  \\
      &  =  &  {\partial}_iA_j-{\partial}_jA_i
      \end{array}   \]
As ${\tilde{C}}_r={\wedge}^rT^*\otimes {\tilde{R}}_2\subset {\wedge}^rT^*\otimes {\hat{R}}_2={\hat{C}}_r$ and ${\hat{R}}_2/{\tilde{R}}_2\simeq T^*$, the conformal Spencer sequence projects onto the sequence $T^*\rightarrow T^*\otimes T^*\rightarrow {\wedge}^2T^*\otimes T^*\rightarrow ...$ which finally projects with a shift by one step onto the Poincare sequence $T^*\stackrel{d}{\rightarrow} {\wedge}^2T^* \stackrel{d}{\rightarrow} {\wedge}^3T^*\rightarrow ... $ by applying the Spencer map $\delta$, because these two sequences are only made by first order involutive operators and are thus formally exact. The short exact sequence $0\rightarrow S_2T^*\stackrel{\delta}{\rightarrow} T^*\otimes T^* \stackrel{\delta}{\rightarrow} {\wedge}^2T^*\rightarrow 0$ has already been used in [23,24] for exhibiting the Ricci tensor and {\it the above result brings for the first time a conformal link between electromagnetism and gravitation} by using second order jets (See [16,17] for more details).  \\
The study of the nonlinear framework is similar. Indeed, using Remark 2.11 with $k=l=r$, we get:\\
\[  {\varphi}_{ij}=A^s_i{\chi}^r_{rs,j}-A^s_j{\chi}^r_{rs,i}={\partial}_i{\chi}^r_{r,j} - {\partial}_j{\chi}^r_{r,i} = {\partial}_i{\alpha}_j-{\partial}_j{\alpha}_i     \]
and we may finish as before as we have taken out the quadratic terms through the contraction.  \\
 \hspace*{12cm}   Q.E.D.  \\ 

This unification result, which may be considered as the ultimate "{\it dream} " of E. and F. Cosserat or H. Weyl, could not have been obtained before 1975 as it can only be produced by means of the (linear/nonlinear) Spencer sequences and NOT by means of the (linear/nonlinear) gauge sequences. \\

\noindent
{\bf 4. Third Part: The Duality Scheme}\\

A duality scheme, first introduced by Henri Poincar\'{e} (1854-1912) in [13], namely a variational framewoirk adapted to the Spencer sequence, could be achieved in local coordinates as we did for the gauge sequence at the end of the introduction. We have indeed presented all the explicit formulas needed for this purpose and the reader will notice that it is difficult or even impossible to find them in [10]. However, it is much more important to relate this dual scheme to homological algebra [25] and algebraic analysis [18,19] by using the comment done at the end of the Second Part which amounts to bring the nonlinear framework to the linear framework, a reason for which the stress equations of continuum mechanics are linear even for nonlinear elasticity [16,22,23].  \\
 
Let $A$ be a {\it unitary ring}, that is $1,a,b\in A \Rightarrow a+b,ab \in A, 1a=a$ and even an {\it integral domain}, that is $ab=0\Rightarrow a=0$ or $b=0$. However, we shall not always assume that $A$ is commutative , that is $ab$ may be different from $ba$ in general for $a,b\in A$. We say that $M={}_AM$ is a {\it left module} over $A$ if $x,y\in M\Rightarrow ax,x+y\in M, \forall a\in A$ or a {\it right module} $M_B$ for $B$ if the operation of $B$ on $M$ is $(x,b)\rightarrow xb, \forall b\in B$. Of course, $A={ }_AA_A$ is a left {\it and} right module over itself. We define the {\it torsion submodule} $t(M)=\{x\in M\mid \exists 0\neq a\in A, ax=0\}\subseteq M$ and $M$ is a {\it torsion module} if $t(M)=M$ or a {\it torsion-free module} if $t(M)=0$. We denote by $hom_A(M,N)$ the set of morphisms $f:M\rightarrow N$ such that $f(ax)=af(x)$. In particular $hom_A(A,M)\simeq M$ because $f(a)=af(1)$ and we recall that a sequence of modules and maps is exact if the kernel of any map is equal to the image of the map preceding it. When $A$ is commutative, $hom(M,N)$ is again an $A$-module for the law $(bf)(x)=f(bx)$ as we have $(bf)(ax)=f(bax)=f(abx)=af(bx)=a(bf)(x)$. In the non-commutative case, things are much more complicate and we have:\\

\noindent
{\bf LEMMA 4.1}: Given ${}_AM$ and ${}_AN_B$, then $hom_A(M,N)$ becomes a right module over $B$ for the law $(fb)(x)=f(x)b$.\\

\noindent
{\it Proof}: We just need to check the two relations:
\[ (fb)(ax)=f(ax)b=af(x)b=a(fb)(x),\]
\[ ((fb')b")(x)=(fb')(x)b"=f(x)b'b"=(fb'b")(x).\]
\hspace*{12cm}               Q.E.D. \\
 
\noindent
{\bf DEFINITION 4.2}: A module $F$ is said to be {\it free} if it is isomorphic to a  power of $A$ called the {\it rank} of $F$ over $A$ and denoted by 
$rk_A(F)$ while the rank of a module is the rank of a maximum free submodule. In the sequel we shall only consider {\it finitely presented} modules, 
namely {\it finitely generated} modules defined by exact sequences of the type $F_1 \stackrel{d_1}{\longrightarrow} F_0 \longrightarrow M\longrightarrow 0$ where $F_0$ and $F_1$ are free modules of finite ranks. For any short exact sequence $0\rightarrow M' \stackrel{f}{\rightarrow} M \stackrel{g}{\rightarrow} M" \rightarrow 0$, we have $rk_A(M)=rk_A(M')+rk_A(M")$. A module $P$ is called {\it projective} if there exists a free module $F$ and another (thus projective) module $Q$ such that $P\oplus Q\simeq F$. A {\it projective (free) resolution} of $M$ is a long exact sequence $... \stackrel{d_3}{\longrightarrow} P_2 \stackrel{d_2}{\longrightarrow} P_1 \stackrel{d_1}{\longrightarrow} P_0 \stackrel{p}{\longrightarrow} M \longrightarrow 0 $ where $P_0, P_1, P_2, ... $ are projective (free) modules, $M=coker(d_1)=P_0/ im(d_1)$ and $p$ is the canonical projection.  \\

We now introduce the {\it extension modules}, using the notation $M^*=hom_A(M,A)$ and, for any morphism $f:M\rightarrow N$, we shall denote by $f^*:N^*\rightarrow M^*$ the morphism which is such that $f^*(h)=h\circ f, \forall h\in hom_A(N,A)$. For this, we take out $M$ in order to obtain the {\it deleted sequence} $... \stackrel{d_2}{\longrightarrow} P_1 \stackrel{d_1}{\longrightarrow} P_0 \longrightarrow 0$ and apply $hom_A(\bullet,A)$ in order to get the sequence $... \stackrel{d^*_2}{\longleftarrow} P^*_1 \stackrel{d^*_1}{\longleftarrow} P^*_0 \longleftarrow 0$. \\

\noindent
{\bf PROPOSITION 4.3}: The extension modules  $ext^0_A(M)=ker(d^*_1)=hom_A(M,A)$ and $ext^i_A(M)=ker(d^*_{i+1})/im(d^*_i), \forall i\geq 1$ do not depend on the resolution chosen and are torsion modules for $i\geq 1$.  \\

Let $\mathbb{Q}\subset K$ be a {\it differential field}, that is a field ($a\in K \Rightarrow 1/a\in K$) with $n$ commuting {\it derivations} $\{{\partial}_1,...,{\partial}_n\}$ with ${\partial}_i{\partial}_j={\partial}_j{\partial}_i={\partial}_{ij}, \forall i,j=1,...,n$ such that ${\partial}_i(a+b)={\partial}_ia+{\partial}_ib$ and ${\partial}_i(ab)=({\partial}_ia)b+a{\partial}_ib, \forall a,b\in K$. Using an implicit summation on multiindices, we may introduce the (noncommutative) {\it ring of differential operators} $D=K[d_1,...,d_n]=K[d]$ with elements $P=a^{\mu}d_{\mu}$ such that $\mu<\infty$ and $d_ia=ad_i+{\partial}_ia$. We notice that $D$ can be generated by $K$ and $T=\{\xi={\xi}^id_i\mid {\xi}^i\in K\}$. Now, if we introduce {\it differential indeterminates} $y=(y^1,...,y^m)$, we may extend $d_iy^k_{\mu}=y^k_{\mu+1_i}$ to ${\Phi}^{\tau}\equiv a^{\tau\mu}_ky^k_{\mu}\stackrel{d_i}{\longrightarrow} d_i{\Phi}^{\tau}\equiv a^{\tau\mu}_ky^k_{\mu+1_i}+{\partial}_ia^{\tau\mu}_ky^k_{\mu}$ for $\tau=1,...,p$. Therefore, setting $Dy^1+...+dy^m=Dy\simeq D^m$, we obtain by residue the {\it differential module} or $D$-{\it module} $M=Dy/D\Phi$. Introducing the two free differential modules $F_0\simeq D^{m_0}, F_1\simeq D^{m_1}$, we obtain equivalently the {\it free presentation} $F_1\stackrel{{\cal{D}}_1}{\rightarrow} F_0 \rightarrow M \rightarrow 0$. More generally, introducing the successive CC as in the preceding section, we may finally obtain the {\it free resolution} of $M$, namely the exact sequence $\hspace{5mm} ... \stackrel{{\cal{D}}_3}{\longrightarrow} F_2  \stackrel{{\cal{D}}_2}{\longrightarrow} F_1 \stackrel{{\cal{D}}_1}{\longrightarrow}F_0\longrightarrow M \longrightarrow 0 $. In actual practice, we let ${\cal{D}}_r$ act on the left on column vectors in the operator case and on the right on row vectors in the module case. Homological algebra has been created for finding intrinsic properties of modules not depending on {\it any} presentation or even on {\it any} resolution.  \\

 We now exhibit another approach by defining the {\it formal adjoint} of an operartor $P$ and an operator matrix $\cal{D}$:  \\

\noindent
{\bf DEFINITION 4.4}: \hspace{1cm}$P=a^{\mu}d_{\mu}\in D  \stackrel{ad}{\longleftrightarrow} ad(P)=(-1)^{\mid\mu\mid}d_{\mu}a^{\mu}   \in D $  \\
\[ <\lambda,{\cal{D}} \xi>=<ad({\cal{D}})\lambda,\xi>+\hspace{1mm} {div}\hspace{1mm} ( ... )  \]
from integration by part, where $\lambda$ is a row vector of test functions and $<  > $ the usual contraction.  \\

\noindent
{\bf LEMMA 4.5}: IIf $f\in aut(X)$, we may set $ x=f^{-1}(y)=g(y)$ and we have the {\it identity}:
\[   \frac{\partial}{\partial y^k}(\frac{1}{\Delta (g(y))} {\partial}_if^k(g(y))\equiv 0.   \]

\noindent
{\bf PROPOSITION 4.6}: If we have an operator $E\stackrel{\cal{D}}{\longrightarrow} F$, we obtain by duality an operator ${\wedge}^nT^*\otimes E^*\stackrel{ad(\cal{D})}{\longleftarrow} {\wedge}^nT^*\otimes F^*$where $E^*$ is obtained from $E$ by inverting the transition matrix. \\

\noindent
{\bf EXAMPLE 4.7}: Let us revisit EM in the light of the preceding results when $n=4$. First of all,  we have $dA=F \Rightarrow dF=0$ in the sequence ${\wedge}^1T^*\stackrel{d}{\longrightarrow} {\wedge}^2T^* \stackrel{d}{\longrightarrow} {\wedge}^3T^*$ and the {\it field equations} are invariant under any local diffeomorphism $f\in aut(X)$. By duality, we get the sequence ${\wedge}^4T^*\otimes {\wedge}^1T \stackrel{ad(d)}{\longleftarrow} {\wedge}^4T^*\otimes {\wedge}^2T \stackrel{ad(d)}{\longleftarrow} {\wedge}^4T^*\otimes {\wedge}^3T$ which is locally isomorphic (up to sign) to ${\wedge}^3T^* \stackrel{d}{\longleftarrow} {\wedge}^2T^* \stackrel{d}{\longleftarrow} {\wedge}^1T^*$ and the {\it induction equations} ${\partial}_i{\cal{F}}^{ij}={\cal{J}}^j$ are thus also invariant under any $f\in aut(X)$. Indeed, using the last lemma and the {\it identity} ${\partial}_{ij}f^l{\cal{F}}^{ij}\equiv 0$, we have: \\
\[\frac{\partial}{\partial y^k}(\frac{1}{\Delta}{\partial}_i f^k{\partial}_j f^l{\cal{F}}^{ij})=\frac{1}{\Delta} {\partial}_i f^k \frac{\partial}{\partial y^k}({\partial}_j f^l{\cal{F}}^{ij})=\frac{1}{\Delta}{\partial}_i({\partial}_j f^l{\cal{F}}^{ij})=\frac{1}{\Delta}{\partial}_j f^l{\partial}_i{\cal{F}}^{ij} \]
Accordingly, it is not correct to say that the conformal group is the biggest group of invariance of Maxwell equations as it is only the biggest group of invariance of the Minkowski constitutive laws in vacuum [12]. Finally, both sets of equations can be parametrized {\it independently}, the first by the potential, the second by the so-called pseudopotential (See [18], p 492 for more details).\\

Now, with operational notations, let us consider the two differential sequences:  \\
\[   \xi  \stackrel{{\cal{D}}}{\longrightarrow} \eta \stackrel{{\cal{D}}_1}{\longrightarrow} \zeta  \]
\[   \nu  \stackrel{ad({\cal{D}})}{\longleftarrow} \mu \stackrel{ad({\cal{D}}_1)}{\longleftarrow} \lambda   \]
where ${\cal{D}}_1$ generates all the CC of ${\cal{D}}$. Then ${\cal{D}}_1\circ {\cal{D}}\equiv 0 \Longleftrightarrow ad({\cal{D}})\circ ad({\cal{D}}_1)\equiv 0 $ but $ad({\cal{D}})$ may not generate all the CC of $ad({\cal{D}}_1)$. Passing to the module framework, we just recognize the definition of $ext^1_D(M)$. Now, exactly like we defined the differential module $M$ from $\cal{D}$, let us define the differential module $N$ from $ad(\cal{D})$. Then $ext^1_D(N)=t(M)$ does not depend on the presentation of $M$ [19]. More generally, changing the presentation of $M$ may change $N$ to $N'$ but we have [11,18]:  \\

\noindent
{\bf THEOREM 4.8}: The modules $N$ and $N'$ are {\it projectively equivalent}, that is one can find two projective modules $P$ and $P'$ such that $N\oplus P\simeq N' \oplus P'$ and we obtain therefore $ext^i_D(N)\simeq ext^i_D(N'), \forall i\geq 1$.  \\

\noindent
{\bf THEOREM 4.9}: When $M$ is a left $D$-module, then $R=hom_K(M,K)$ is also a left $D$-module.  \\

\noindent
{\it Proof}:  Let us define:\\
\[   (af)(m)=af(m)=f(am) \hspace{1cm} \forall a\in K, \forall m\in M\]
\[   (\xi f)(m)=\xi f(m)-f(\xi m)  \hspace{1cm}  \forall \xi ={\xi}^id_i\in T, \forall m\in M  \]
It is easy to check that $d_ia=ad_i+{\partial}_ia$ in the operator sense and that $\xi\eta -\eta\xi =[\xi,\eta]$ is the standard bracket of vector fields. We finally get $(d_if)^k_{\mu}=(d_if)(y^k_{\mu})={\partial}_if^k_{\mu}-f^k_{\mu +1_i}$ that is {\it exactly} the {\it Spencer operator} we used in the second part. In fact, $R$ is the {\it projective limit} of ${\pi}^{q+r}_q:R_{q+r}\rightarrow R_q$ in a coherent way with jet theory [18,19].  \\
\hspace*{12cm}     Q.E.D.   \\

\noindent
{\bf COROLLARY 4.10}: if $M$ and $N$ are right $D$-modules, then $hom_K(M,N)$ becomes a left $D$-module.  \\

\noindent
{\it Proof}: We just need to set $(\xi f)(m)=f(m\xi)-f(m)\xi, \forall \xi\in T, \forall m\in M$ and conclude as before.  \\
\hspace*{12cm}  Q.E.D.  \\

As $D={ }_DD_D$ is a bimodule, then $M^*=hom_D(M,D)$is a right $D$-module according to Lemma 4.1 and the module $N_r$ defined by the ker/coker sequence $0\longleftarrow N_r \longleftarrow F^*_1 \stackrel{{\cal{D}}^*}{\longleftarrow} F^*_0 \longleftarrow M^* \longleftarrow 0$ is in fact a right module $N_r=N_D$.\\

\noindent
{\bf THEOREM 4.11}: We have the {\it side changing} procedure $N=N_l={ }_DN=hom_K({\wedge}^nT^*,N_r)$.  \\

\noindent
{\it Proof}: According to the above Corollary, we just need to prove that ${\wedge}^nT^*$ has a natural right module structure over $D$. For this, if $\alpha=adx^1\wedge ...\wedge dx^n\in T^*$ is a volume form with coefficient $a\in K$, we may set $\alpha.P=ad(P)(a)dx^1\wedge...\wedge dx^n$ when $P\in D$. As $D$ is generated by $K$ and $T$, we just need to check that the above formula has an intrinsic meaning for any $\xi={\xi}^id_i\in T$. In that case, we check at once:
\[  \alpha.\xi=-{\partial}_i(a{\xi}^i)dx^1\wedge...\wedge dx^n=-\cal{L}(\xi)\alpha \]
by introducing the Lie derivative of $\alpha$ with respect to $\xi$, along the intrinsic formula ${\cal{L}}(\xi)=i(\xi)d+di(\xi)$ where $i( )$ is the interior multiplication and $d$ is the exterior derivative of exterior forms. According to well known properties of the Lie derivative, we get :
\[\alpha.(a\xi)=(\alpha.\xi).a-\alpha.\xi(a), \hspace{5mm} \alpha.(\xi\eta-\eta\xi)=-[\cal{L}(\xi),\cal{L}(\eta)]\alpha=-\cal{L}([\xi,\eta])\alpha=\alpha.[\xi,\eta].  \]
\hspace*{12cm}  Q.E.D.  \\
 
\noindent
{\bf REMARK 4.12}: The above results provide a new light on duality in physics. Indeed, as the Poincar\'{e} sequence is self-adjoint (up to sign) {\it as a whole} and the linear Spencer sequence for a Lie group of transformations is locally isomorphic to copies of that sequence, it follows from Proposition 4.3 that $ad(D_{r+1})$ parametrizes $ad(D_r)$ in the dual of the linear Spencer sequence while $ad({\cal{D}}_{r+1})$ parametrizes $ad({\cal{D}}_r)$ in the dual of the linear Janet sequence, a result highly not evident at first sight in view of the Janet/Spencer diagram for the conformal group of tranformations of space-time that we have presented because ${\cal{D}}_r$ and $D_{r+1}$ are {\it totally different operators}.  \\
 
\noindent
{\bf 5. Conclusion}  \\
 
The mathematical foundations of Gauge Theory (GT) leading to Yang-Mills equations are {\it always} presented in textbooks or papers without quoting that the group theoretical methods involved are exactly the same as the standard ones used in continuum mechanics, particularly in the analytical mechanics of rigid bodies and in hydrodynamics. Surprisingly, the lagrangians of GT are (quadratic) functions of the curvature $2$-form while the lagrangians of mechanics are (quadratic or cubic) functions of the potential $1$-form. Meanwhile, the corresponding variational principle leading to Euler-Lagrange equations is also {\it shifted by one step in the use of the same gauge sequence}. This situation is contradicting the well known field/matter couplings existing between elasticity and electromagnetism (piezzoelectricity, photoelasticity). In this paper, we prove that the mathematical foundations of GT are not coherent with jet theory and the Spencer sequence. Accordingly, they must be revisited within this new framework, that is when there is {\it a Lie group of transformations considered as a Lie pseudogroup}, contrary to the situation existing in GT. Such a new approach, based on new mathematical tools still not known today by physicists, allows to unify electromagnetism and gravitation. Finally, the striking fact that the Cosserat/Maxwell/Weyl equations can be parametrized, contrary to Einstein equations, is shown to have quite deep roots in homological algebra through the use of {\it extension modules} and {\it duality theory} in the framework of {\it algebraic analysis}.  \\

\noindent
{\bf REFERENCES}  \\

\noindent
[1] V. ARNOLD: M\'{e}thodes Math\'{e}matiques de la M\'{e}canique Classique, Appendice 2 (G\'{e}od\'{e}siques des m\'{e}triques invariantes \`{a} gauche sur des groupes de Lie et hydrodynamique des fluides parfaits), MIR, moscow, 1974,1976. \\
\noindent
[2] V. ARNOLD: Sur la G\'{e}ometrie des Groupes de Lie de Dimension Infini et ses Applications \{a} l'Hydrodynamique des Fluides Parfaits, Ann. Inst. Fourier, Grenoble, 16 (1),1966, pp. 319-361.  \\
\noindent
[3]Ê. BIRKHOFF: Hydrodynamics, Princeton University Press, 1954.  \\
\noindent
[4] D. BLEECKER: Gauge Theory and Variational Principles, Addison-Wesley, 1981; Dover, 2005.  \\
\noindent
[5] E. COSSERAT, F. COSSERAT: Th\'{e}orie des Corps D\'{e}formables, Hermann, Paris, 1909.\\
\noindent
[6] W. DRECHSLER,  M.E. MAYER: Fiber Bundle Techniques in Gauge Theories, Springer Lecture Notes in Physics 67, Springer, 1977.  \\
\noindent
[7] M. GOCKELER: Differential Geometry, Gauge Theories and Gravity, Cambridge Monographs on Mathematical Physics, Cambride University Press, 1987.  \\
\noindent
[8] M. JANET: Sur les Syst\`{e}mes aux D\'{e}riv\'{e}es Partielles, Journal de Math., 8, 1920, pp. 65-151. \\
\noindent
 [9] S. KOBAYASHI, K. NOMIZU: Foundations of Differential Geometry, Vol I, J. Wiley, New York, 1963, 1969.\\
\noindent
[10] A. KUMPERA, D.C. SPENCER: Lie Equations, Ann. Math. Studies 73, Princeton University Press, Princeton, 1972.\\
\noindent
[11] E. KUNZ: Introduction to Commutative Algebra and Algebraic Geometry, BirkhaŸser, 1985.  \\
\noindent
[12] V. OUGAROV: Th\'{e}orie de la Relativit\'{e} Restreinte, MIR, Moscow, 1969 ( french, 1979).\\
\noindent
[13] H. POINCARE: Sur une Forme Nouvelle des Equations de la M\'{e}canique, C. R. Acad\'{e}mie des Sciences Paris, 132 (7), 1901, pp. 369-371.  \\
\noindent
[14] J.-F. POMMARET: Systems of Partial Differential Equations and Lie Pseudogroups, Gordon and Breach, New York, 1978; Russian translation: MIR, Moscow, 1983.\\
\noindent
[15] J.-F. POMMARET: Differential Galois Theory, Gordon and Breach, New York, 1983.\\
\noindent
[16] J.-F. POMMARET: Lie Pseudogroups and Mechanics, Gordon and Breach, New York, 1988.\\
\noindent
[17] J.-F. POMMARET: Partial Differential Equations and Group Theory, Kluwer, 1994.\\
http://dx.doi.org/10.1007/978-94-017-2539-2    \\
\noindent
[18] J.-F. POMMARET: Partial Differential Control Theory, Kluwer, Dordrecht, 2001.\\
\noindent
[19] J.-F. POMMARET: {\it Algebraic Analysis of Control Systems Defined by Partial Differential Equations}, Advanced Topics in Control Systems Theory, Springer, Lecture Notes in Control and Information Sciences 311, 2005, Chapter 5, pp. 155-223.\\
\noindent
[20] J.-F. POMMARET: Group Interpretation of Coupling Phenomena, Acta Mechanica, 149, 2001, pp. 23-39.\\
http://dx.doi.org/10.1007/BF01261661  \\
\noindent
[21] J.-F. POMMARET: Arnold's Hydrodynamics Revisited, AJSE-mathŽmatiques, 1, 1, 2009, pp. 157-174.  \\
\noindent
[22] J.-F. POMMARET: Parametrization of Cosserat Equations, Acta Mechanica, 215, 2010, pp. 43-55.\\
http://dx.doi.org/10.1007/s00707-010-0292-y  \\
\noindent
[23] J.-F. POMMARET: Spencer Operator and Applications: From Continuum Mechanics to Mathematical Physics, in "Continuum Mechanics-Progress in Fundamentals and Engineering Applications", Dr. Yong Gan (Ed.), ISBN: 978-953-51-0447--6, InTech, 2012, Available from: \\
http://www.intechopen.com/books/continuum-mechanics-progress-in-fundamentals-and-engineering-applications/spencer-operator-and-applications-from-continuum-mechanics-to-mathematical-physics  \\
\noindent
[24] J.-F. POMMARET: The Mathematical Foundations of General Relativity Revisited, Journal of Modern Physics, 2013, 4, pp. 223-239. \\
 http://dx.doi.org/10.4236/jmp.2013.48A022   \\
\noindent
[25] J. J. ROTMAN: An Introduction to Homological Algebra, Academic Press, 1979.\\
\noindent
[26] D. C. SPENCER: Overdetermined Systems of Partial Differential Equations, Bull. Am. Math. Soc., 75, 1965, pp. 1-114.\\
\noindent
[27] P.P. TEODORESCU: Dynamics of Linear Elastic Bodies, Editura Academiei, Bucuresti, Romania; Abacus Press, Tunbridge, Wells, 1975.\\
\noindent
[28] E. VESSIOT: Sur la Th\'{e}orie des Groupes Infinis, Ann. Ec. Norm. Sup., 20, 1903, pp. 411-451.\\
\noindent
[29] H. WEYL: Space, Tilme, Matter, Springer, 1918, 1958; Dover, 1952. \\
\noindent
[30] C.N. YANG: Magnetic Monopoles, Fiber Bundles and Gauge Fields, Ann. New York Acad. Sciences, 294, 1977, pp. 86-97.  \\
\noindent
[31] C.N. YANG: R.L. MILLS: Conservation of Isotopic Gauge Invariance, Phys. Rev., 96, 1954, pp. 191-195.\\
\noindent
[32] Z. ZOU, P. HUANG, Y. ZHANG, G. LI: Some Researches on Gauge Theories of Gravitation, Scientia Sinica, XXII, 6, 1979, pp. 628-636.\\

\end{document}